\documentclass{article}
\usepackage[utf8]{inputenc}

\title{One Hundred Probability and Statistics Inequalities}
\author{CNP Slagle}
\date{\today}

\usepackage{Definitions}
\begin{document}

\maketitle

\setlength{\parindent}{0cm}

\newtheorem{Lemma}{Lemma}
\newtheorem{Theorem}{Theorem}
\maketitle
\section{Introduction}
In 2012, the author compiled a subset of the following inequalities for a researcher in randomized algorithms.  One might think of said inequalities as a very quick reference, with access to primary and secondary resources listed either within the section or alongside the inequality of interest.  In the intervening years, some of the original sources and their respective links have vanished, leading the author to consider a companion document with proofs for select inequalities within this list.  Though the author would refute the completeness of this collection for more advanced researchers, he nonetheless believes it may serve some interest.  It is important  
\section{Basic Probability and Measure Theory Inequalities}
The relations to follow include axioms within the probabilistic framework, along with a few of the basic inferences derived therefrom.  See \cite{CasellaBerger:1} for a wonderful introduction.
Given events (sets) $A$, $B$, and countable $\{A_{n}\}_{n=1}^{\infty}$,
\begin{enumerate}
\item $\PP[A]\geq 0$
\item $\PP[A]\leq 1$
\item If $A\subset B$, then $\PP[A]\leq \PP[B]$
\item If $A\subset B$, then $\PP[B]\leq \PP[A^{c}]$
\item \textbf{(Boole)} $\PP[\cup_{n=1}^{\infty}A_{n}]\leq \sum_{n=1}^{\infty}\PP[A_{n}]$
\item $\PP[\cup_{n=1}^{\infty}A_{n}]\geq \sup\{\PP[A_{n}]|n=1,2,\dots\}$
\item $\PP[\cap_{n=1}^{\infty}A_{n}]\leq \inf\{\PP[A_{n}]|n=1,2,\dots\}$
\item $\PP[A \cap B] \leq \min\{\PP[A],\PP[B]\}$
\item \textbf{(Bonferroni)} $\PP[A \cap B] \geq \PP[A]+\PP[B]-1$ 
\item \textbf{(Bonferroni General)} $\PP[\cap_{i=1}^{n}] \geq \sum_{i=1}^{n}\PP[A_{i}]-(n-1)$
\item $\PP[A|B] \geq \PP[A\cap B]$
\item \textbf{(Karlin Ost)} Define $P_{1}=\sum_{i=1}^{n}\PP[A_{i}]$, $P_{2}=\sum_{1\leq i < j\leq n}\PP[A_{i}\cap A_{j}]$, $P_{3}=\sum_{1\leq i < j < k\leq n}\PP[A_{i}\cap A_{j} \cap A_{k}]$,...,$P_{n}=\PP[A_{1}\cap\dots\cap A_{n}]$.  Then for $i=1,\dots,n$,
\begin{equation}
P_{1}-P_{2}+P_{3}-\dots\pm P_{i-1}\geq \PP\left[\cup_{i=1}^{n}A_{i}\right] \geq P_{1}-P_{2}+P_{3}-\dots\mp P_{i}.
\end{equation}

\end{enumerate}
\section{Power Means \cite{CasellaBerger:1}}
Define the \textit{pth power mean} of a finite set of positive numbers $S$ to be
\begin{equation}
PM(S,p)=\sqrt[p]{\sum_{s\in S}{\frac{{s}^p}{|S|}}}
\end{equation}
Notice that the arithmetic mean and harmonic mean of the set $S$ are simply $A_{M}(S)=PM(S,1)$ and $H_{M}(S)=PM(S,-1)$, respectively.  Less clear is that the geometric mean $G_{M}(S) = \sqrt[|S|]{\prod_{s\in S}s}=\limit{p\to 0}PM(S,p)$, the maximum $\max(S)=\limit{p\to \infty}PM(S,p)$, and the minimum $\min(S)=\limit{p\to -\infty}PM(S,p)$.
So, we have
\begin{enumerate}
\item $PM(S,p_{0}) \leq PM(S,p_{1})$ for $p_{0}\leq p_{1}$ 
\item $H_{M}(S)\leq G_{M}(S)\leq A_{M}(S)$ 
\end{enumerate}

\section{Expectations and Variances I}
The following inequalities range from elementary to moderate complexity, all available in \cite{CasellaBerger:1}.  Let $X$ and $Y$ be random variables.  If an inequality includes a function $f$ of a random variable $X$, assume that the expectation $\EE f(X)$ exists.  
\begin{enumerate}
\item If $g(X)\leq h(X)$, then $\EE g(X) \leq \EE h(X)$. 
\item If $a \leq g(X) \leq b$, then $a \leq \EE g(X) \leq b$.
\item \textbf{(H\"{o}lder)} If $p,q$ satisfy ${{1}\over{p}}+{{1}\over{q}}=1$, then $|\EE XY|\leq \EE|XY|\leq (\EE|X|^{p})^{\frac{1}{p}}(\EE|Y|^{q})^{{1}\over{q}}$
\item \textbf{(Jensen)} For a convex function $g$, If $X\geq Y$, then $\EE g(X)\geq g(EX)$. 
\item \textbf{(Cauchy-Schwartz)} $|\EE XY|\leq \EE|XY|\leq \sqrt{(\EE|X|^{2})(\EE|Y|^{2})}$ 
\item $Var(X)\geq 0$
\item $Cov^{2}(X,Y)\leq Var(X)Var(Y)$
\item \textbf{(H\"{o}lder Special Case)} For $p>1$, $\EE |X|\leq \sqrt[p]{\EE|X|^{p}}$
\item \textbf{(Liapounov)} For $s>r>1$, $\sqrt[r]{\EE|X|^{r}}\leq\sqrt[s]{\EE|X|^{s}}$
\item \textbf{(Minkowski)} For $p \geq 1$, $\sqrt[p]{\EE|X+Y|^{p}}\leq \sqrt[p]{E|X|^{p}}+\sqrt[p]{E|Y|^{p}}$
\item \textbf{(Triangle)} As a special case of Minkowski's inequality, $\EE |X+Y|\leq E|X| + E|Y|$.
\item If $g$ is nondecreasing and $h$ is nonincreasing, then $\EE (g(X)h(X))\leq(Eg(X))(Eh(X))$.
\item If $g$ and $h$ are both nondecreasing or both nonincreasing, then $\EE (g(X)h(X))\geq(Eg(X))(Eh(X))$.

\item \textbf{(Cram\'{e}r-Rao)} Suppose $X_{1},\dots,X_{n}$ is a sample with joint pdf $f(\mathbf{x}|\theta)$ and $W(\mathbf{X})$ is any estimator of $\theta$ such that ${{d}\over{d\theta}}E_{\theta}W(\mathbf{X})=\int_{\chi}{{\partial}\over{\partial\theta}}[W(\mathbf{x})]f(\mathbf{x}|\theta)d\mathbf{x}$ and $Var_{\theta}(W(\mathbf{X}))<\infty$.  Then
\begin{equation}
Var_{\theta}(W(\mathbf{X}))\geq {{({{d}\over{d\theta}}E_{\theta}W(\mathbf{X}))^{2}}\over{E_{\theta}\left[\left({{\partial}\over{\partial\theta}}\log(f(\mathbf{x}|\theta))\right)^{2}\right]}}.
\end{equation}

\item \textbf{(Cram\'{e}r-Rao IID)} Suppose $X_{1},\dots,X_{n}$ is a sample iid with marginal pdf $f(x|\theta)$ and $W(\mathbf{X})$ is any estimator of $\theta$ such that ${{d}\over{d\theta}}E_{\theta}W(\mathbf{X})=\int_{\chi}{{\partial}\over{\partial\theta}}[W(\mathbf{x})]f(\mathbf{x}|\theta)d\mathbf{x}$ and $Var_{\theta}(W(\mathbf{X}))<\infty$.  Then
\begin{equation}
Var_{\theta}(W(\mathbf{X}))\geq {{({{d}\over{d\theta}}E_{\theta}W(\mathbf{X}))^{2}}\over{nE_{\theta}\left[\left({{\partial}\over{\partial\theta}}\log(f(x|\theta))\right)^{2}\right]}}.
\end{equation}

\item \textbf{(Rao-Blackwell)} Let $U$ be an unbiased estimator of $\tau(\theta)$, and let $T$ be a sufficient statistic for $\theta$.  Define $\phi(T)=E(U|T)$.  Then $\EE \phi(T)=\tau(\theta)$, and
\begin{equation}
Var_{\theta}\phi(T)\leq Var_{\theta}W\text{ for all }\theta.
\end{equation}
\end{enumerate}

\section{Expectations and Variances II}
\begin{enumerate}
\item \textbf{(Han \cite{BoucheronBousquetLugosi:1})} Let $X_{1},\dots,X_{n}$ be independent discrete random variables. Let $H(X_{\pi_{1}},\dots,X_{\pi_{k}})$ be the joint entropy of a subset of the $\{X_{i}\}$.  Then
\begin{equation}
H(X_{1},\dots,X_{n})\leq \sum_{i=1}^{n} H(X_{1},\dots,X_{i-1},X_{i+1},\dots,X_{n}).  
\end{equation}

\item \cite{BoucheronBousquetLugosi:1} Let $X_{1},\dots,X_{n}$ be independent random variables.  Let $g:Domain(X_{1},\dots,X_{n})\to\mathbb{R}$ be Lesbegue measurable, and $Z=g(X_{1},\dots,X_{n})$. Then
\begin{equation}
Var(Z)\leq \sum_{i=1}^{n} E[(Z-E(Z|X_{1},\dots,X_{i-1},X_{i+1},\dots,X_{n})^{2}].
\end{equation}

\item \textbf{(Efron-Stein \cite{BoucheronBousquetLugosi:1})} Let $X_{1},\dots,X_{n}$ be independent random variables.  Let $g:Domain(X_{1},\dots,X_{n})\to\mathbb{R}$ be Lesbegue measurable, and $Z=g(X_{1},\dots,X_{n})$. Let $Y_{1},\dots,Y_{n}$ be an independent copy of $X_{1},\dots,X_{n}$, and let $Z_{i}=g(X_{1},\dots,Y_{i},\dots,X_{n})$.  Then
\begin{equation}
Var(Z)\leq \sum_{i=1}^{n} E[(Z-Z_{i})^{2}]. 
\end{equation}
\item \textbf{(Logarithmic Sobolev \cite{BoucheronBousquetLugosi:1})} Let $X_{1},\dots,X_{n}$ be independent random variables.  Let \\$g_{i}:Domain(X_{1},\dots,X_{i-1},X_{i+1},\dots,X_{n})\to\mathbb{R}$ be Lesbegue measurable, $Z_{i}=g_{i}(X_{1},\dots,X_{i-1},X_{i+1},\dots,X_{n})$, $g:Domain(X_{1},\dots,X_{n})\to\mathbb{R}$ be Lesbegue measurable, and $Z=g(X_{1},\dots,X_{n})$. Let $\psi(t)=e^{t}-t-1$ and $s>0$.  Then
\begin{equation}
sE(Ze^{sZ})-E(e^{sZ})\log[E(e^{sZ})]\leq \sum_{i=1}^{n} E[e^{sZ}\psi(-s(Z-Z_{i}))].  
\end{equation}
\item \textbf{(Symmetrized Logarithmic Sobolev \cite{BoucheronBousquetLugosi:1})} Let $X_{1},\dots,X_{n}$ be independent random variables.  Let \\$g:Domain(X_{1},\dots,X_{n})\to\mathbb{R}$ be Lesbegue measurable, and $Z=g(X_{1},\dots,X_{n})$. Let $Y_{1},\dots,Y_{n}$ be an independent copy of $X_{1},\dots,X_{n}$, and let $Z_{i}=g(X_{1},\dots,Y_{i},\dots,X_{n})$.  Let $\psi(t)=e^{t}-t-1$ and $s>0$.  Then
\begin{equation}
sE(Ze^{sZ})-E(e^{sZ})\log[E(e^{sZ})]]\leq \sum_{i=1}^{n} E[e^{sZ}\psi(-s(Z-Z_{i}))]. 
\end{equation}
\item \cite{Wasserman:2} Suppose $\{X_{n}\}$ is a sequence of random variables such that for all $n$, $X_{n}\geq0$, and for all $\epsilon>0$, there exist $c_{1}>0$ and $c_{2}>e^{-1}$ such that $\PP[X_{n}>\epsilon]\leq c_{1}e^{-c_{2}n\epsilon^{2}}$.  Then
\begin{equation}
\EE X_{n}\leq \sqrt{{1+\log(c_{1})}\over{nc_{2}}}.
\end{equation} 

\item \textbf{(Ledoux-Talagrand Contraction \cite{Duchi:1})} Suppose $X_{i},\dots,X_{n}$ are iid \textit{Rademacher} variables ($\PP[X_{i}=1]=\PP[X_{i}=-1]={{1}\over{2}}$).  Suppose $f:\mathbb{R}^{+}\to\mathbb{R}^{+}$ be convex and increasing, and $\phi_{i}:\mathbb{R}\to\mathbb{R}$ be Lipschitz with constant $L$ for $i=1,\dots,n$.  Then for $T\subset\mathbb{R}^{n}$,
\begin{equation}
Ef\left({{1}\over{2}}\sup_{\mathbf{t}\in T}\left|\sum_{i=1}^{n}X_{i}\phi_{i}(t_{i})\right|\right) \leq Ef\left(L \sup_{\mathbf{t}\in T}\left|\sum_{i=1}^{n}X_{i}t_{i}\right|\right) 
\end{equation}
\item \textbf{(Bhatia-Davis \cite{Bhatia:1})} If a univariate probability distribution $F$ has minimum $m$, maximum $M$, and mean $\mu$, then for any $X$ following $F$, $Var(X) \leq (M-\mu)(\mu-m)$. 

\item \textbf{(Popoviciu \cite{Vasile:1})} If a univariate probability distribution $F$ has minimum $m$ and maximum $M$, then for any $X$ following $F$, $Var(X) \leq  {1\over4}(M-m)^{2}$. 

\item \textbf{(Chapman-Robbins \cite{Chapman:1})} Suppose $\mathbf{X}$ is a random variable in $\mathbb{R}^{k}$ with an unknown parameter $\theta$.  If $\delta(\mathbf{X})$ is an unbiased estimator for $\tau(\theta)$, then
\begin{equation}
Var(\delta(\mathbf{X}))\geq \sup_{\Delta} {{[\tau(\theta+\Delta)-\tau(\theta)]^{2}}\over{E_{\theta}\left[{{p(\mathbf{X}.\theta+\Delta)}\over{p(\mathbf{X},\theta)}}-1\right]^{2}}}.
\end{equation}

\item \textbf{(Entropy Power \cite{Dembo:1})} Define the entropy of $X$ to be $\entropy(X)=-\EE\log f_{X}(X)$, where $f_{X}(x)$ is the pdf or pmf of $X$.  Define the entropy power of $X$ to be $N(X)={{1}\over{2\pi e}}e^{{2\over n}h(X)}$.  Then for random variables $X$ and $Y$, we have $N(X+Y)\geq N(X)+N(Y)$.

\item \textbf{(Marcinkiewicz Zygmund \cite{Marcinkiewicz:1})} Let $X_{1},\dots,X_{n}$ be independent random variables with common support such that $\EE X_{i}=0$ and $\EE X_{i}^{p}<\infty$ for all $p\geq1$.  Then there exist constants $A(p)$ and $B(p)$, dependent only on $p$, such that
\begin{equation}
A(p)E\left(\sum_{i=1}^{n}|X_{i}|^2\right)^{p\over 2} \leq E\left(\sum_{i=1}^{n}|X_{i}|\right)^{p}\leq B(p)E\left(\sum_{i=1}^{n}|X_{i}|^2\right)^{p\over 2}.  \end{equation}
\item \textbf{(Khintchine \cite{Wolff:1})} Let $X_{1},\dots,X_{n}$ be iid Rademacher random variables.  Then for any $\lambda_{1},\dots,\lambda_{n}\in \mathbb{C}$ and $p>0$, there exist constants $A(p)$ and $B(p)$, dependent only on $p$, such that
\begin{equation}
A(p)\left(\sum_{i=1}^{n}|\lambda_{i}|^2\right)^{1\over 2} \leq \left(E\left(\sum_{i=1}^{n}|\lambda_{i}X_{i}|\right)^{p}\right)^{1\over p}\leq B(p)\left(\sum_{i=1}^{n}|\lambda_{i}|^2\right)^{1\over 2}.
\end{equation}

\item \textbf{(Rosenthal I \cite{Rosenthal:1})} Let $X_{1},\dots,X_{n}$ be independent nonnegative random variables such that $\EE X_{i}^{p}<\infty$ for a fixed $p\geq1$, $i=1,\dots,n$. Then there exist constants $A(p)$ and $B(p)$ dependent only on $p$ such that
\begin{equation}
\begin{array}{rl}
A(p)\max\left\{\sum_{i=1}^{n}EX_{i}^{p},\left(\sum_{i=1}^{n}EX_{i}\right)^{p}\right\}\leq & \left(\sum_{i=1}^{n}EX_{i}\right)^{p}\\\\
\leq  & B(p)\max\left\{\sum_{i=1}^{n}EX_{i}^{p},\left(\sum_{i=1}^{n}EX_{i}\right)^{p}\right\}.
\end{array}
\end{equation}

\item \textbf{(Rosenthal II)} Let $X_{1},\dots,X_{n}$ be independent random variables such that $\EE X_{i}=0$ and $\EE X_{i}^{p}<\infty$ for a fixed $p\geq1$, $i=1,\dots,n$. Then there exist constants $A(p)$ and $B(p)$ dependent only on $p$ such that
\begin{equation}
\begin{array}{rl}
A(p)\max\left\{\sum_{i=1}^{n}E|X_{i}|^{p},\left(\sum_{i=1}^{n}EX_{i}^{2}\right)^{p\over 2}\right\} & \leq \left|\sum_{i=1}^{n}EX_{i}\right|^{p}\\\\
 & \leq B(p)\max\left\{\sum_{i=1}^{n}E|X_{i}|^{p},\left(\sum_{i=1}^{n}EX_{i}^{2}\right)^{p\over 2}\right\}.
\end{array}
\end{equation}

\item \textbf{(Papadatos \cite{Papadatos:1})} Let $X_{(1)},\dots,X_{(n)}$ be the order statistics of iid random variables $X_{1},\dots,X_{n}$ with variance $\sigma^{2}$.  Define $G(x)=I_{x}(k,n+1-k)$ and $\sigma_{n}^{2}(k)=\sup_{0<x<1}\left[{{G(x)(1-G(x))}\over{x(1-x)}}\right]$.  Then
\begin{equation}
Var(X_{(k)})\leq \sigma_{n}^{2}(k)\sigma^{2}.
\end{equation}

\item \textbf{(H\"{u}rlimann Upper \textit{n--r} \cite{Hurlimann:1})} Let $X_{(1)},\dots,X_{(n)}$ be the order statistics of iid random variables $X_{1},\dots,X_{n}$.  Define $\bar{X}={1\over n}\sum_{i=1}^{n}X_{i}$, and the \textit{biased observed variance} $S^{2}={1\over n}\sum_{i=1}^{n}(X_{i}-\bar{X})^{2}$.  Then for $r=0,\dots,n-1$, the \textit{average of the upper n-r order statistics} satisfies
\begin{equation}
{1\over{n-r}}\sum_{i=r+1}^{n}X_{(i)}\leq \bar{X}+S\sqrt{{r}\over{n-r}}.\end{equation}

\item \textbf{(H\"{u}rlimann Average Excess \cite{Hurlimann:1})} Let $X_{(1)},\dots,X_{(n)}$ be the order statistics of iid random variables $X_{1},\dots,X_{n}$.  Define $\bar{X}={1\over n}\sum_{i=1}^{n}X_{i}$, and the \textit{biased observed variance} $S^{2}={1\over n}\sum_{i=1}^{n}(X_{i}-\bar{X})^{2}$.  Then for $r=0,\dots,n-1$, the \textit{average excess of the upper $n-r$ order statistics conditioned on the $r$th order statistic} satisfies
\begin{equation}
{1\over{n-r}}\sum_{i=r+1}^{n}(X_{(i)}-X_{(r)})\leq S{{n}\over{\sqrt{r(n-r)}}}.
\end{equation}
\item \textbf{(H\"{u}rlimann Stop-Loss Excess \cite{Hurlimann:1})} Let $X_{(1)},\dots,X_{(n)}$ be the order statistics of iid random variables $X_{1},\dots,X_{n}$.  Define $\bar{X}={1\over n}\sum_{i=1}^{n}X_{i}$, and the \textit{biased observed variance} $S^{2}={1\over n}\sum_{i=1}^{n}(X_{i}-\bar{X})^{2}$.  Define the \textit{rth stop-loss statistics} to be $SL(d_{r})=\sum_{i=1}^{n}(X_{(i)}-(n-r)d_{r})$ for $d_{r}\in [X_{(r)},X_{(r+1)}]$.  Then for $r=0,\dots,n-1$,
\begin{equation}
SL(d_{r})\leq (n-r)\left[\bar{X}-d_{r}+S\sqrt{r\over{n-r}}\right].  
\end{equation}
\end{enumerate}

\section{Concentration Inequalities}
We may apply concentration inequalities to bound probabilities, often for sequences of random variables.  Such inequalities bound tail probabilities, though the applicability and accuracy of such depends highly upon the local and asymptotic behavior of the random variables of interest.  Let $X$ be a random variable.
\subsection{Elementary to Moderate Concentration}
\begin{enumerate}
\item (Markov\cite{CasellaBerger:1}) Suppose $X\ge 0$, and $\EE[X]>0$.  Then $\PP(X\ge t)\le \frac{\EE[X]}{t}$ for all $t>0$.

\item \textbf{(Chebychev \cite{CasellaBerger:1})} For $t>0$, $\PP[|X-EX|\geq t]\leq {{Var(X)}\over{t^2}}$.

\item \textbf{($g$-Markov)} Let $X\ge 0$, $\EE[X]>0$.  Then for increasing $g:\RR^+_0\to\RR^+$, 
\begin{equation}
    \PP(X\ge t)\le \frac{\EE[g(X)]}{g(t)}.
\end{equation}

\item \textbf{(Normal I, Mill \cite{CasellaBerger:1}\cite{Wasserman:1})} For $Z$ a standard normal, $\PP[|Z|\geq t]\leq \sqrt{{2}\over{\pi}}{e^{-t^{2}/2}\over{t}}$.

\item \textbf{(Normal II \cite{CasellaBerger:1})} For $Z$ a standard normal, $\PP[|Z|\geq t]\geq \sqrt{{2}\over{\pi}}e^{-t^{2}/2}{t\over{1+t^2}}$.

\item \textbf{(Chernoff I \cite{CasellaBerger:1}\cite{Mitzenmacher:1})} Let $M_{X}(t)$,$-h\leq t\leq h$ be the moment-generating function of $X$.  Then $\PP[X> a]\leq e^{-at}M_{X}(t)$, $-h\leq t \leq h$.

\item \textbf{(Chernoff II \cite{CasellaBerger:1}\cite{Mitzenmacher:1})} Let $M_{X}(t)$,$-h\leq t\leq h$ be the moment-generating function of $X$.  Then $\PP[X\leq a]\leq e^{-at}M_{X}(t)$, $-h\leq t \leq 0$.

\item \textbf{(Chernoff Sum I \cite{CasellaBerger:1})} Let $X_{1},\dots,X_{n}$ be iid, $X=\sum_{i=1}^{n}X_{i}$, and $M_{X}(t)$,$-h\leq t\leq h$ be the moment-generating function of $X_1$.  .  Then $\PP[S>a]\leq e^{-at}[M_{X}(t)]^{n}$ for $0\leq t\leq h$.

\item \textbf{(Chernoff Sum II \cite{CasellaBerger:1})} Let $X_{1},\dots,X_{n}$ be iid, $X=\sum_{i=1}^{n}X_{i}$, and $M_{X}(t)$,$-h\leq t\leq h$ be the moment-generating function of $X_1$.  Then $\PP[S\leq a]\leq e^{-at}[M_{X}(t)]^{n}$ for $-h\leq t\leq 0$.

\item \textbf{(Chernoff Mean \cite{CasellaBerger:1})} Let $X_{1},\dots,X_{n}$ be iid, $\epsilon>0$, $\bar{X_{n}}=\sum_{i=1}^{n}X_{i}$, $M_{U}(t)$, $-h_{U}\leq t\leq h_{U}$ be the moment-generating function of $U=X_{1}-EX_{1}-\epsilon$, and $M_{V}(t)$, $-h_{V}\leq t\leq h_{V}$ be the moment-generating function of $V=-X_{1}+EX_{1}-\epsilon$.  Then there exist for some $0 < t_{U}\leq h_{U}$ and $-h_{V}\leq t_{V}<0$  \footnote{Such a $t_{U}$ and $t_{V}$ exist since $\EE U < 0$ and $\EE V<0$, guaranteeing that $M_{U}$ and $M_{V}$ are decreasing in a neighborhood of zero.} such that
\begin{equation}
\PP[|\bar{X_{n}}-EX_{1}|>\epsilon]\leq 2c^{n},\\ 
\text{ where }c=\max\{M_{U}(t_{U}),M_{V}(t_{V})\} \in (0,1) . 
\end{equation}
\end{enumerate}

\subsection{Randomized Algorithms}
See \cite{Mitzenmacher:1} for an introduction into randomized algorithms, whence we infer the following inequalities.
\begin{enumerate}
\item \textbf{(Chernoff Poisson Trials I)} Let $X_{i}$ be $n$ independent Poisson trials  \footnote{Each $X_{i}$ is a Bernoulli$(p_{i})$.}.  Let $X=\sum_{i=1}^{n}X_{i}$.  Then for $\delta>0$,
\begin{equation}
\PP[X\geq (1+\delta)EX] < \left({{e^{\delta}}\over{(1+\delta)^{1+\delta}}}\right)^{EX}.
\end{equation}

\item \textbf{(Chernoff Poisson Trials II)} Let $X_{i}$ be $n$ independent Poisson trials.  Let $X=\sum_{i=1}^{n}X_{i}$.  Then for $0<\delta\leq 1$,
\begin{equation}
\PP[X\geq (1+\delta)EX] < e^{-(EX)\delta^{2}/3}.
\end{equation}

\item \textbf{(Chernoff Poisson Trials III)} Let $X_{i}$ be $n$ independent Poisson trials.  Let $X=\sum_{i=1}^{n}X_{i}$.  Then for $R\geq gEX$,
\begin{equation}
\PP[X\geq R] < 2^{-R}.
\end{equation}

\item \textbf{(Chernoff Poisson Trials IV)} Let $X_{i}$ be $n$ independent Poisson trials.  Let $X=\sum_{i=1}^{n}X_{i}$.  Then for $0<\delta<1$,
\begin{equation}
\PP[X\leq (1-\delta)EX] < \left({{e^{-\delta}}\over{(1-\delta)^{1-\delta}}}\right)^{EX}.
\end{equation}
\item \textbf{(Chernoff Poisson Trials V)} Let $X_{i}$ be $n$ independent Poisson trials.  Let $X=\sum_{i=1}^{n}X_{i}$.  Then for $0<\delta<1$,
\begin{equation}
\PP[X\leq (1-\delta)EX] < e^{-\delta^{2}EX/2}.
\end{equation}

\item \textbf{(Chernoff Rademacher I)}  Suppose $X_{1},\dots,X_{n}$ be iid such that $\PP[X_{i}=1]=\PP[X_{i}=-1]={{1}\over{2}}$.  If $X=\sum_{i=1}^{n}X_{i}$ and $a>0$, then $\PP[X\geq a]\leq e^{{-a^{2}}\over{2n}}$.

\item \textbf{(Chernoff Rademacher II)}  Suppose $X_{1},\dots,X_{n}$ be iid such that $\PP[X_{i}=1]=\PP[X_{i}=-1]={{1}\over{2}}$.  If $X=\sum_{i=1}^{n}X_{i}$ and $a>0$, then $\PP[|X|\geq a]\leq 2e^{{-a^{2}}\over{2n}}$.

\item \textbf{(Chernoff Bernoulli I)}  Suppose $X_{1},\dots,X_{n}$ be iid Bernoulli$\left({1\over2}\right)$. If $X=\sum_{i=1}^{n}X_{i}$ and $0<a<{n\over2}$, then $\PP\left[X\leq {n\over2} - a\right]\leq 2e^{{-2a^{2}}\over{n}}$.

\item \textbf{(Chernoff Bernoulli II)}  Suppose $X_{1},\dots,X_{n}$ be iid Bernoulli$\left({1\over2}\right)$. If $X=\sum_{i=1}^{n}X_{i}$ and $0<\delta<1$, then $\PP\left[X\leq {n\over2}(1-\delta)\right]\leq 2e^{{-n\delta^{2}}\over{2}}$. 

\item \textbf{(Chernoff Bernoulli III)}  Suppose $X_{1},\dots,X_{n}$ be iid Bernoulli$\left({1\over2}\right)$. If $X=\sum_{i=1}^{n}X_{i}$ and $a>0$, then $\PP\left[X\geq {n\over2} + a\right]\leq 2e^{{-2a^{2}}\over{n}}$.

\item \textbf{(Chernoff Bernoulli IV)}  Suppose $X_{1},\dots,X_{n}$ be iid Bernoulli$\left({1\over2}\right)$. If $X=\sum_{i=1}^{n}X_{i}$ and $\delta>0$, then $\PP\left[X\leq {n\over2}(1+\delta)\right]\leq 2e^{{-n\delta^{2}}\over{2}}$.

\end{enumerate}

\subsection{Unimodality, and Misc.}
\begin{enumerate}
\item \textbf{(Gauss \cite{CasellaBerger:1})} Suppose $X$ follows a unimodal distribution with mode $\nu$, and define $\tau^{2}=E(X-\nu)^{2}$. Then
\begin{equation}
\PP[|X-\nu|>\epsilon]\leq 
\begin{cases}
{{4\tau^2}\over{9\epsilon^2}}, & \epsilon\geq \sqrt{{4}\over{3}}\tau \\
1-{{\epsilon}\over{\tau \sqrt{3}}}, & \epsilon\leq \sqrt{{4}\over{3}}\tau 
\end{cases}
\end{equation}

\item \textbf{(Vysochanski\u{i}-Petunin \cite{CasellaBerger:1})} Suppose $X$ follows a unimodal distribution, and define $\xi^{2}=E(X-\alpha)^{2}$ for arbitrary $\alpha$.  Then 
\begin{equation}
\PP\left[|X-\alpha|>\epsilon\right]\leq 
\begin{cases}
{{4\xi^2}\over{9\epsilon^2}}, & \epsilon\geq \sqrt{{8}\over{3}}\xi \\
{{4\xi^2}\over{9\epsilon^2}}-{{1}\over{3}}, & \epsilon\leq \sqrt{{8}\over{3}}\xi 
\end{cases}
\end{equation}

\item \textbf{(Hoeffding I \cite{Wasserman:1})} Let $Y_{1},\dots,Y_{n}$ be independent observations such that $\EE Y_{i}=0$ and $a_{i}\leq Y_{i} \leq b_{i}$ for all $i$.  If $\epsilon >0$ and $t>0$, then
\begin{equation}
\PP\left[\sum_{i=1}^{n}Y_{i}\geq\epsilon \right]\leq e^{-t\epsilon}\prod_{i=1}^{n}e^{t^{2}(b_{i}-a_{i})^{2}/8}
\end{equation}

\item \textbf{(Hoeffding II \cite{Wasserman:1})} Let $X_{1},\dots,X_{n}$ be independent Bernoulli($p$).  If $\epsilon >0$, then
\begin{equation}
\PP\left[\left|\sum_{i=1}^{n}X_{i}-np\right|\geq\epsilon \right]\leq 2e^{-2n\epsilon^2}
\end{equation}

\item \textbf{(Saw)} Suppose $X_{1},\dots,X_{n}$ are iid with finite first and second order moments.  Let $\bar{X}={{1}\over{n}}\sum_{i=1}^{n}X_{i}$ and $S^{2}={{1}\over{n-1}}\sum_{i=1}^{n}(X_{i}-\bar{X})^{2}$.  Let $k>0$, $\nu(t)=\max\left\{m\in\mathbb{N}|m<{{n+1}\over{t}}\right\}$, $\alpha(t)={{(n+1)(n+1-\nu(t))}\over{1+\nu(t)(n+1-\nu(t))}}$, and $\beta={{n(n+1)k^{2}}\over{n-1+(n+1)k^{2}}}$.  Then
\begin{equation}
\PP[|X-\bar{X}|\geq kS]\leq
\begin{cases}
{{1}\over{n+1}}(\nu(\beta)-1) & \text{ if }\nu\text{ is odd and }\beta>\alpha(\beta) \\
{{1}\over{n+1}}\nu(\beta) & \text{ otherwise.}
\end{cases}
\end{equation}

\end{enumerate}

\section{Kannan Combinatoric Inequalities}
Kannan \cite{Kannan:1} furnishes an array of inequalities helpful in analyzing graphs and other objects of combinatoric import.

\begin{enumerate}
    \item \textbf{(Chromatic Number)} Let $G(\{1,\dots,n\},P)$ be a random graph with edge probabilities $P=\left(p_{ij}\right)$.  The \textit{chromatic number} $\chi=\chi(G)$ is the least number of colors necessary to color $G$ such that no two vertices sharing an edge receive the same color.  Let $p={{\sum_{i,j}p_{ij}}\over{n\choose 2}}$.  Then there exists a constant $c>0$ such that for $t\in(0,n\sqrt{p})$,
\begin{equation}
\PP[|\chi(G)-E\chi(G)|\geq t]\leq e^{{-ct^{2}}\over{n\sqrt{p}\log n}}.
\end{equation} 

\item \textbf{(Johnson-Lindenstrauss Random Projection)} Suppose $k\leq n$, and we pick $V_{1},\dots,V_{k}$ uniformly randomly from the surface of the unit ball in $\mathbb{R}^{n}$.  Then for $\epsilon\in(0,1)$, there exist constants $c_{1},c_{2}>0$ such that
\begin{equation}
\PP\left[\left|\sum_{i=1}^{k}v_{i}^{2}-{{k}\over{n}}\right|\geq {{\epsilon k}\over{n}}\right]\leq c_{1}e^{-c_{2}\epsilon^{2}k}
\end{equation}

\item \textbf{(Random Projection)} Suppose $m$ is an even positive integer and $X_{1},\dots,X_{n}$ are real-valued random observations satisfying the \textit{strong negative correlation principle}.  That is, for all $i$, $\EE X_{i}(X_{1}+\dots+X_{i-1})^{l}<0$ when $l<m$ is odd and $\EE (X_{i}^{l}|X_{1}+\dots+X_{i-1})\leq\left({{n}\over{m}}\right)^{{l-2}\over{2}}l!$ for $l\leq m$ even.  Define constants $\{M_{i,l}\}$, $\{K_{i,l}\}$, and $\{L_{i,l}\}$ such that $\EE (X_{i}^{l}|X_{1}+\dots+X_{i-1})\leq M_{i,l}$, each $K_{i,l}$ is an indicator variable on the \textit{typical} case of the conditional expectation where $\PP[K_{i,l}]=1-\delta_{i,l}$, and $\EE (X_{i}^{l}|X_{1}+\dots+X_{i-1},K_{i,l})\leq L_{i,l}$ for $l=2,4,\dots,m$ and $i=1,\dots,n$.  Finally, let $X=\sum_{i=1}^{n}X_{i}$. Then
\begin{equation}
\begin{array}{rl}
    \EE X^{m}\leq & (cm)^{{{m}\over{2}} + 1}\left(\sum_{l=1}^{m/2}{{m^{1-{{1}\over{l}}}\over{l^{2}}}}\left(\sum_{i=1}^{n}L_{i,2l}\right)^{1\over l}\right)^{m\over 2}\\\\
  + & (cm)^{m+2}\sum_{l=1}^{m/2}{{1}\over{nl^{2}}}\sum_{i=1}^{n}\left(nM_{i,2l}\delta_{i,2l}^{2/(m-2l+2)}\right)^{2\over{ml}}.
\end{array}
 \end{equation}
\item \textbf{(Bin Packing)} Suppose $Y_{1},\dots,Y_{n}$ are iid from a discrete distribution of $r$ atoms each with probability at least $1\over \log n$ and $\EE Y_{1} \leq {{1}\over{r^{2}\log n}}$.  Let $f(Y_{1},\dots,Y_{n})$ be the minimum number of unit capacity bins necessary to pack the $Y_{1},\dots,Y_{n}$ items.  Then there exist constants $c_{1},c_{2}>0$ such that if $t\in(0,n[(EY_{i})^{3}+Var(Y_{i})])$, then
\begin{equation}
\PP[|f-Ef|\geq t+r] \leq c_{1}e^{{-c_{2}t^{2}}\over{n[(EY_{i})^{3}+Var(Y_{i})]}}. \end{equation}

\item \textbf{(Strong Negative Correlation)} Suppose $m$ is an even positive integer, and $X_{1},\dots,X_{n}$ are real-valued random observations satisfying the \textit{strong negative correlation principle}.  That is, for all $i$, $\EE X_{i}(X_{1}+\dots+X_{i-1})^{l}<0$ when $l<m$ is odd and $\EE (X_{i}^{l}|X_{1}+\dots+X_{i-1})\leq\left({{n}\over{m}}\right)^{{l-2}\over{2}}l!$ for $l\leq m$ even.  Then
\begin{equation}
\EE\left(\sum_{i=1}^{n}X_{i}\right)^{m}\leq (24mn)^{{m}\over{2}}.
\end{equation}
\item \textbf{(Hamiltonian Tour)} Suppose $Y_{1},\dots,Y_{n}$ are sets of points generated independently and respectively from $n$ subsquares of size ${{1}\over{\sqrt{n}}}\times{1\over\sqrt{n}}$ of the unit square, and there exists a constant $c\in(0,1)$ such that $\PP[|Y_{i}|]\leq c$ for all $i$.  Suppose further that for $\epsilon>0$, and $l\in\{1,\dots,{m\over 2}\}$, $\EE |Y_{i}|^{l}\leq [O(l)]^{(2-\epsilon)l}$.  Finally, suppose $f(Y_{1},\dots,Y_{n})$ is the length of the shortest Hamiltonian tour through $Y_{1}\cup\dots\cup Y_{n}$.  Then
\begin{equation}
\EE[f(Y_{1},\dots,Y_{n})-Ef(Y_{1}\dots,Y_{n})]^{m}\leq (cm)^{{m}\over{2}}.
\end{equation}
\item \textbf{(MST)} Suppose $Y_{1},\dots,Y_{n}$ are sets of points generated independently and respectively from $n$ subsquares of size ${{1}\over{\sqrt{n}}}\times{1\over\sqrt{n}}$ of the unit square, and there exists a constant $c\in(0,1)$ such that $\PP[|Y_{i}|]\leq c$ for all $i$.  Suppose further that for $\epsilon>0$, and $l\in\{1,\dots,{m\over 2}\}$, $\EE |Y_{i}|^{l}\leq [O(l)]^{(2-\epsilon)l}$.  Finally, suppose $f(Y_{1},\dots,Y_{n})$ is the length of a minimum spanning tree of $Y_{1}\cup\dots\cup Y_{n}$.  Then
\begin{equation}
E[f(Y_{1},\dots,Y_{n})-Ef(Y_{1}\dots,Y_{n})]^{m}\leq (cm)^{{m}\over{2}}.
\end{equation}
\item \textbf{(Random Vector)} Suppose $\mathbf{Y}=(Y_{1},\dots,Y_{n})$ is a random vector such that for a fixed $k\leq n$, $\EE (Y_{i}^{2}|Y_{1}^{2},\dots,Y_{i}^{2})$ is a nondecreasing function of $Y_{1}^{2}+\dots+Y_{i-1}^2$ for $i=1,\dots,k$ and for even $l\leq k$, there exists a $c>0$ such that $\EE (Y_{i}^{l}|Y_{1}^{2},\dots,Y_{i-1}^{2})\leq \left({{cl}\over{n}}\right)^{l\over 2}$.  Then for any even $m\leq k$,
\begin{equation}
\EE\left(\sum_{i=1}^{k}Y_{i}^{2}-EY_{i}^{2}\right)^{m}\leq \left({{\sqrt{cmk}}\over{n}}\right)^{m}. 
\end{equation}

\end{enumerate}

\section{Means and Variances II}
The inequalities to follow furnish mechanisms for the analysis of interdependence, Markov chains, vectors, and graphs, among others.
\begin{enumerate}
\item \textbf{(Talagrand \cite{Tao:1})} Let $\mathbf{X}$ be chosen randomly uniformly from $\{-1,1\}^{n}$, let $A$ be a convex subset of $\mathbb{R}^{n}$, $A_{t}=\{\mathbf{p}\in\mathbb{R}^{n}|dist(\mathbf{p},A)\leq t\}$.  Then there exists $c>0$ such that $\PP[\mathbf{X}\in A]\PP[\mathbf{X}\notin A_{t}]\leq e^{-ct^{2}}$ for all $t>0$.

\item \textbf{(Talagrand Large Deviation \cite{Tao:1})} Let $\mathbf{X}$ be chosen randomly uniformly from $\{-1,1\}^{n}$, $V$ be a $d$-dimensional subspace of $\mathbb{R}^{n}$.  Then there exist constants $c,C>0$ such that $\PP[|dist(X,V)-\sqrt{n-d}|\geq t]\leq Ce^{-ct^{2}}$ for all $t>0$.

\item \textbf{(Gaussian for Lipschitz \cite{Tao:1})} Let $\mathbf{X}$ be an $n$-dimensional random vector such that each $X_{i}$ is an independent $n(0,1)$ variable.  If $f:\mathbb{R}^{d}\to \mathbb{R}$ is a Lipschitz function with scale constant 1\footnote{A Lipschitz function $f$ satisfies $|f(x)-f(y)|\leq M||x-y||$ for all $x,y\in domain(f)$.}, then there exists a constant $c>0$ such that $\PP[|f(\mathbf{X})-Ef(\mathbf{X})| \geq t] \leq e^{-ct^{2}}$ for all $t>0$.

\item \textbf{(Azuma \cite{Chung:1})} Suppose $X_{0},\dots,X_{n}$ is a \textit{martingale} ($\EE (X_{i}|X_{1},\dots,X_{i-1})=X_{i-1}$ for $i=1,\dots,n$); suppose further that $X$ is $\mathbf{c}$-Lipschitz ($|X_{i}-X_{i-1}|\leq c_{i}$ for $i=1,\dots,n$, $c\in\mathbb{R}^{n}$ positive); then
\begin{equation}
\PP[X_{n}-X_{0}\geq\lambda]\leq 2e^{{-\lambda^{2}}\over{2\sum_{i=1}^{n}c_{i}^2}}.
\end{equation} 

\item \textbf{(Bennett \cite{BoucheronBousquetLugosi:1})} Let $X_{1},\dots,X_{n}$ be independent random variables of zero mean such that $\PP[X_{i}\leq 1]=1$.  Let $h(u)=(1+u)\log(1+u)-u$ for $u\geq0$ and $\sigma^{2}={{1}\over{n}}\sum_{i=1}^{n}Var(X_{i})$.  Then for $t>0$,
\begin{equation}
\PP\left[\sum_{i=1}^{n}X_{i}>t\right]\leq e^{-n\sigma^{2}h\left({{t}\over{n\sigma^{2}}}\right)}. \end{equation}

\item \textbf{(Bernstein \cite{BoucheronBousquetLugosi:1})} Let $X_{1},\dots,X_{n}$ be independent random variables of zero mean such that $\PP[X_{i}\leq 1]=1$.  Let $\sigma^{2}={{1}\over{n}}\sum_{i=1}^{n}Var(X_{i})$.  Then for $\epsilon>0$,
\begin{equation}
\PP\left[{{1}\over{n}}\sum_{i=1}^{n}X_{i}>\epsilon\right]\leq e^{{-n\epsilon^{2}}\over{2(\sigma^{2}+\epsilon/3)}}.
\end{equation}

\item \textbf{(McDiarmid Bounded Differences I \cite{Bartlett:1})} Let $X_{1},\dots,X_{n}$ be independent random variables each whose domain is $\chi$.  If $f:\chi^{n}\to\mathbb{R}^{n}$ is a function such that for all $\mathbf{x}\in\chi^{n}$, $y\in\chi$, and $i\in\{1,\dots,n\}$, there exists a constant $c_{i}>0$ such that $|f(\mathbf{x})-f(x_{1},\dots,x_{i-1},y,x_{i+1},\dots,x_{n})|\leq c_{i}$, then
\begin{equation}
\PP[f(\mathbf{X})-Ef(\mathbf{X})\geq t] \geq  e^{{-t^{2}}\over{\sum_{i=1}^{n}c_{i}^{2}}} \text{ for all $t>0$}.
\end{equation}

\item \textbf{(McDiarmid Bounded Differences II \cite{Bartlett:1})} Let $X_{1},\dots,X_{n}$ be independent random variables each whose domain is $\chi$.  If $f:\chi^{n}\to\mathbb{R}^{n}$ is a function such that for all $\mathbf{x}\in\chi^{n}$, $y\in\chi$, and $i\in\{1,\dots,n\}$, there exists a constant $c_{i}>0$ such that $|f(\mathbf{x})-f(x_{1},\dots,x_{i-1},y,x_{i+1},\dots,x_{n})|\leq c_{i}$, then
\begin{equation}
\PP[f(\mathbf{X})-Ef(\mathbf{X})\leq -t] \geq  e^{{-t^{2}}\over{\sum_{i=1}^{n}c_{i}^{2}}} \text{ for all $t>0$.}
\end{equation}

\item \textbf{(Dvoretzky Kiefer Wolfowitz I \cite{DvoretzkyKieferWolfowitz:1})} Suppose $X_{1},\dots,X_{n}$ are iid univariate random variables following cdf $F$.  Let $F_{n}(x)={1\over n}\sum_{i=1}^{n}1_{X_{i}\leq x}$ be the empirical distribution.  Then for $\epsilon>\sqrt{{1\over 2n} \log 2}$,
\begin{equation}
\PP\left[\sup_{x\in \mathbb{R}}(F_{n}(x)-F(x))>\epsilon\right]\leq e^{-2n\epsilon^{2}}.
\end{equation}

\item \textbf{(Dvoretzky Kiefer Wolfowitz II \cite{DvoretzkyKieferWolfowitz:1})} Suppose $X_{1},\dots,X_{n}$ are iid univariate random variables following cdf $F$.  Let $F_{n}(x)={1\over n}\sum_{i=1}^{n}1_{X_{i}\leq x}$ be the empirical distribution.  Then for $\epsilon>0$,
\begin{equation}
\PP\left[\sup_{x\in \mathbb{R}}|F_{n}(x)-F(x)|>\epsilon\right]\leq 2e^{-2n\epsilon^{2}}.
\end{equation}

\item \textbf{(Etemadi Differing Means \cite{Etemadi:1})} Let $X_{1},\dots,X_{n}$ be random variables with common support.  Let $S_{k}=\sum_{i=1}^{k}X_{k}$ be the $k$th partial sum.  Then for $\epsilon>0$,
\begin{equation}
\PP\left[\max_{1\leq k\leq n}|S_{k}|\geq 3\epsilon\right]\leq 3\max_{1\leq k\leq n}\PP\left[|S_{k}|\geq \epsilon\right].
\end{equation}

\item \textbf{(Etemadi Shared Means \cite{Etemadi:1})} Let $X_{1},\dots,X_{n}$ be random variables with common support and equal means.  Let $S_{k}=\sum_{i=1}^{k}X_{k}$ be the $k$th partial sum.  Then for $\epsilon>0$,
\begin{equation}
\PP\left[\max_{1\leq k\leq n}|S_{k}|\geq \epsilon\right]\leq {27\over\epsilon^{2}}Var(S_{n}).
\end{equation}

\item \textbf{(Kolmogorov \cite{Billingsley:1})} Let $X_{1},\dots,X_{n}$ be independent random variables with common support such that $\EE X_{i}=0$ and $Var(X_{i})<\infty$ for $i=1,\dots,n$.  Let $S_{k}=\sum_{i=1}^{k}X_{k}$ be the $k$th partial sum.  Then for $\epsilon>0$,
\begin{equation}
\PP\left[\max_{1\leq k\leq n}|S_{k}|\geq \epsilon\right]\leq {1\over\epsilon^{2}}\sum_{i=1}^{n}Var(X_{i}). \end{equation}

\item \textbf{(Chebychev Multidimensional \cite{Vershynin:1})} Let $\mathbf{X}\in\mathbb{R}^{n}$ be a random vector with covariance matrix \\$V=E\left[(\mathbf{X}-E\mathbf{X})(\mathbf{X}-E\mathbf{X})^{T}\right]$.  Then for $t>0$,
\begin{equation}
\PP\left[\sqrt{(\mathbf{X}-E\mathbf{X})^{T}V^{-1}(\mathbf{X}-E\mathbf{X})}\right]\leq {n\over t^{2}}.
\end{equation}

\item \textbf{(Leguerre Samuelson \cite{Samuelson:1})}  Let $X_{1},\dots,X_{n}$ be random variables with common support, and define $\bar{X}={1\over n}\sum_{i=1}^{n}X_{i}$ and $S^{2}={1\over{n-1}}\sum_{i=1}^{n}(X_{i}-\bar{X})^{2}$.  Then for $i=1,\dots,n$ with probability one,
\begin{equation}
\bar{X}-S\sqrt{n-1}\leq X_{i} \leq \bar{X}+S\sqrt{n-1}.
\end{equation}

\item \textbf{(LeCam \cite{LeCam:1})} Suppose $X_{1},\dots,X_{n}$ are independent binomial random variables with respective success parameters $p_{1},\dots,p_{n}$.  Letting $\lambda_{n}=\sum_{i=1}^{n}p_{i}$, we have
\begin{equation}
\sum_{k=0}^{\infty}\left|\PP\left[\sum_{i=1}^{n}X_{i}=k\right]-{\lambda_{n}^{k}e^{-\lambda_{n}}\over{k!}}\right|\leq 2\sum_{i=1}^{n}p_{i}^2.
\end{equation}

\item \textbf{(Doob Martingale \cite{RevuzYor:1})} Let $\mathbf{X}\in\mathbb{R}^{n}$ be a \textit{martingale} ($\EE (X_{i}|X_{1},\dots,X_{i-1})=X_{i-1}$ for $i=2,\dots,n$).  Then for $C>0$, $p\geq 1$,
\begin{equation}
\PP\left[\sup_{1\leq i \leq n}X_{i}\geq C\right]\leq {{EX_{n}^{p}}\over{C^{p}}}.
\end{equation}
\end{enumerate}

\bibliography{inequalities}
\bibliographystyle{ieeetr}

\end{document}